\documentstyle{article}
\title{The transversal wave equation and the noncommutative geodesic flow
in Riemannian foliations}
\author{Yuri A. Kordyukov\\Department of Mathematics,\\
Ufa State Aviation Technical University,\\K.Marx str. 12,
Ufa 450000, Russia\\ e-mail: yurikor@math.ugatu.ac.ru}
\date{January 28, 1997}
\def\Bbb#1{{\bf #1}}
\newtheorem{thm}{Theorem}[section]

\newtheorem{cor}[thm]{Corollary}
\newtheorem{defn}[thm]{Definition}

\newtheorem{lem}[thm]{Lemma}

\newtheorem{prop}[thm]{Proposition}

\newtheorem{rem}[thm]{Remark}

\newtheorem{ack}{Acknowledgment}

\newcommand{\be}{\begin{equation}}
\newcommand{\ee}{\end{equation}}
\newcommand{\bea}{\begin{eqnarray}}
\newcommand{\eea}{\end{eqnarray}}
\newcommand{\ba}{\begin{eqnarray*}}
\newcommand{\ea}{\end{eqnarray*}}
\newcommand{\bl}{\begin{lem}}
\newcommand{\el}{\end{lem}}
\newcommand{\bt}{\begin{thm}}
\newcommand{\et}{\end{thm}}
\newcommand{\bc}{\begin{cor}}
\newcommand{\ec}{\end{cor}}
\newcommand{\bp}{\begin{prop}}
\newcommand{\ep}{\end{prop}}
\newcommand{\p}{{\em Proof}. }

\newcommand{\F}{{\cal F}}

\newcommand{\supp}{\mbox{\rm supp}\;}

\newcommand{\cinf}{C^\infty}
\begin{document}
\maketitle
\section{Introduction}
\label{Intr}
Let $(M,{\cal F})$ be a compact Riemannian foliated manifold, equipped
with a bundle-like metric $g_M$.  We will also use the following
notation: $N^*{\cal F}$ is the conormal bundle to ${\cal F}$,
$G$ is the holonomy groupoid of ${\cal F}$.

In Introduction, we will formulate our main results for the
geometric case of the transverse signature operator, referring the
reader to the main body of the paper for the formulations in the
case of a general transversally elliptic operator.
Recall that the transverse signature operator is a first order differential
operator in the space $C^{\infty}(M,\Lambda T^{*}M)$:
$$
D_H=d_H + d^*_H,
$$
where $d_H$ is the $(0,1)$-component of the de Rham differential $d$
in the bigrading of $\Lambda T^{*}M$:
$$\Lambda^k T^{*}M=\bigoplus_{i+j=k}
\Lambda^{i}T{\F}^{*}\otimes\Lambda^{j}(T{\F}^{\bot})^{*}.
$$
We will consider $D_H$ as an operator, acting in the space 
$C^{\infty}(M,\Lambda^{j}(T{\F}^{\bot})^{*})$ of transversal differential 
forms. 

Let  $g_t$ be the geodesic flow
of the Riemannian metric $g_M$, and, for any $t\in {\Bbb R}$, 
$\Lambda(t)$ be
the graph of the symplectomorphism $g_t$ in $\tilde{T}^*M=T^*M\setminus\{0\}$:
$$
\Lambda(t)=\{((x,\xi),(y,\eta))\in 
\tilde{T}^*M\times \tilde{T}^*M : (x,\xi)=g_t(y,\eta)\}.
$$
\bt
\label{wavelapl}
The operator $e^{it|D_H|}$ has the form
$$
e^{it|D_H|}=W(t)+R(t),
$$
where $W(t)$ is a Fourier integral operator associated with the canonical
relation $\Lambda(t)$, $W(t)\in I^0(M\times M, 
\Lambda'(t);\Omega^{1/2}(M\times M))$,
$R(t)$ is a smooth family of operators from 
${\rm OP}^{2,-\infty}(M,N^*{\cal F},E)$ (see Definition~\ref{op}).
\et

The first application of Theorem~\ref{wavelapl} concerns to the singularities
of the Fourier transform of the distributional spectrum
distribution function of the operator $|D_H|$. 
Since the bundle $\Lambda^*(T{\F}^{\bot})^{*}$ is holonomy 
equivariant, there is defined
a $\ast$-representation $R$ of the involutive
algebra $C^{\infty}_c(G)$ in $L^2(M,\Lambda^{j}(T{\F}^{\bot})^{*})$ 
(see (\ref{repr})).
Let a function $\theta(t)\in {\cal D}'(G), t\in {\Bbb R}$ be given 
by the formula
$$
<\theta(t),k>= {\rm tr}\; R(k)e^{it|D_H|}, k\in C^{\infty}_c(G).
$$
Since the metric $g_M$ on $M$ is bundle-like, the geodesic flow 
preserves $N^*{\cal F}$, inducing a flow on $N^*{\cal F}$, denoted also
by $g_t$.
 
The following theorem extends the result on
a relationship between the spectrum of Laplacian on a compact manifold
and the set of lengths of closed geodesics due to Chazarain,
Duistermaat-Guillemin \cite{Chaz,DG} to the case of transversally elliptic
operators on a Riemannian foliated manifold. The similar result for 
transversally globally elliptic operators in ${\Bbb R}^n$ is stated 
in \cite{HR}. 

\bt
\label{DG}
For any $k\in C^{\infty}_c(G)$, the function $\theta_k(t)=<\theta(t),k>$
is smooth
outside the set of $t\in {\Bbb R}$ such that there exist $\gamma
\in \supp k, \nu\in N^*_{r(\gamma)}{\cal F}\setminus \{0\}$, with
$$
g_{-t}\ dh_{\gamma}^*(\nu)=\nu,
$$
where $dh_{\gamma}: N^*_{s(\gamma)}{\cal F}\rightarrow 
N^*_{r(\gamma)}{\cal F}$ is the linear holonomy map, corresponding to 
$\gamma$.
\et
Otherwise speaking, the singularities of $\theta_k$ are contained in the 
set of lengths of transversal (that is, orthogonal to leaves) 
geodesics $l:[0,1]\rightarrow M$ from 
$x$ to $y$ such that $x$ and 
$y$ are in the same leaf and there exists a  $\gamma\in {\rm supp}\;k, 
\gamma : x\rightarrow y$ such that the corresponding linear
holonomy map maps the velocity vector $\dot{l}(0)$ to $l$ at $x$ 
to the velocity vector $\dot{l}(1)$ to $l$ at $y$. 
Since $k$ is compactly supported, this set is discrete 
in ${\Bbb R}$.   

The second application of Theorem~\ref{wavelapl} is a Egorov type 
theorem for transversally elliptic operators. 
The theorem concerns to a description of the operator 
$$
\Phi_t(K)=e^{it|D_H|}Ke^{-it|D_H|},
$$ 
where $K$ is a transversal pseudodifferential operator of 
class $\Psi^{m,-\infty}(M,{\cal F},E)$ introduced in \cite{noncom}.

The corresponding symbol classes are defined as follows.
Let ${\cal F}_N$ be a foliation in $\tilde{N}^*{\cal F}$,
which is the lift of ${\cal F}$
via a leafwise flat (Bott) connection in $\tilde{N}^*{\cal F}$ \cite{Molino}.
Let $G_{{\cal F}_N}$ be the holonomy groupoid  of the lifted
foliation ${\cal F}_N$.
The symbol class $S^{m}(G_{{\cal F}_N},{\cal L}(\pi^*E)\otimes\Omega^{1/2})$
is the space of all
$s\in C^{\infty}(G_{{\cal F}_N},{\cal L}(\pi^*E)\otimes\Omega^{1/2})$
homogeneous of degree $m$ such that $\pi_G(\supp s)$ is compact
in $G_{\cal F}$ (here $\Omega^{1/2}=r^*(\Omega^{1/2}({\cal F}_N))\otimes
(s^*(\Omega^{1/2}({\cal F}_N)))^*$, with $\Omega^{1/2}({\cal F}_N)$, being 
the vector bundle of leafwise half-densities on $N^*{\cal F}$). The space
$S^{0}(G_{{\cal F}_N},{\cal L}(\pi^*E)\otimes\Omega^{1/2})$ carries
the structure of an involutive algebra, which provides a description
of the transversal unitary cotangent space $SN^*{\cal F}/{\cal F}_N$
from the point of noncommutative geometry.
The algebra $\Psi^{*,-\infty}(M,{\cal F},E)$ of transversal
pseudodifferential operators and the algebra
$S^{m}(G_{{\cal F}_N},{\cal L}(\pi^*E)\otimes\Omega^{1/2})$ of
transversal symbols are connected by
the principal symbol mapping
\be
\label{symb}
\Psi^{m,-\infty}(M,{\cal F},E)\rightarrow
S^{m}(G_{{\cal F}_N},{\cal L}(\pi^*E)\otimes\Omega^{1/2}).
\ee
Let $g_t$ denote also the restriction of the geodesic flow to $N^*{\F}$.
There is a flow $\{G_t\}$ on $G_{{\cal F}_N}$ such that
$ G_t\circ s=G_t\circ r=g_t$.
The transversal geodesic flow  is  
the one-parameter group $G^*_t$ of automorphisms of the
involutive algebra 
$S^m(G_{{\cal F}_N},{\cal L}(\pi^*E)\otimes\Omega^{1/2})$,
induced by $G_t$. From the point of noncommutative geometry,
$G^*_t$ is an analogue of the (classical) geodesic flow for the
transversal unitary cotangent space $\tilde{N}^*{\cal F}/{\cal F}_N$.

\bt
\label{Egorov0}
Given an operator $K\in \Psi^{m,-\infty}(M,{\cal F},E)$ with the principal
symbol $k\in S^{m}(G_{{\cal F}_N},{\cal L}(\pi^*E)\otimes\Omega^{1/2})$,
the operator 
$\Phi_t(K)$ 
is an operator of class
$\Psi^{m,-\infty}(M,{\cal F},E)$ with the principal symbol
$$
k(t)=G^*_t(k)\in S^m(G_{{\cal F}_N},{\cal L}(\pi^*E)\otimes\Omega^{1/2}).
$$
\et

Theorem~\ref{Egorov0} is applied to describe the noncommutative
geodesic flow in transversal geometry of foliations.
Consider a spectral triple
$({\cal A},{\cal H},D)$ (in a sense of \cite{Co-M,Sp-view}) studied in 
\cite{noncom}:
the involutive
algebra ${\cal A}$ is the algebra $\cinf_c(G)$, the Hilbert space ${\cal H}$
is the $L^2$-space $L^2(M,\Lambda^*(T{\F}^{\bot})^*)$ of transversal
differential forms on which
an element $k$ of the algebra ${\cal A}$ is represented via
the $\ast$-representation $R$, and the operator $D$ is the 
tranverse signature operator $D_H$.

Alain Connes \cite{Sp-view} defined the notions of the unitary 
cotangent bundle $S^*{\cal A}$,
being an involutive algebra, and of the noncommutative geodesic flow
$\alpha_t$, being a one-parameter group of automorphisms of $S^*{\cal A}$
 (see Section~\ref{ncg} for details).
In the case under consideration, we prove existence of an
surjective projection of
$S^*{\cal A}$ to 
$\bar{S}^{0}(G_{{\cal F}_N},{\cal L}(\pi^*E)\otimes\Omega^{1/2})$
(given by the principal symbol map (\ref{symb}) in transversal $\Psi$DO
calculus):
$$
P: S^*{\cal A}\rightarrow
\bar{S}^{0}(G_{{\cal F}_N},{\cal L}(\pi^*E)\otimes\Omega^{1/2}),
$$
and Theorem~\ref{Egorov0} states that, under this projection,
the noncommutative geodesic flow $\alpha_t$,
defined by the triple $({\cal A}, {\cal H}, D)$, maps to the transversal
geodesic flow $G^*_t$. More precisely, we have the following theorem:
\bt
\label{noncom:geo}
Given the spectral triple $({\cal A},{\cal H},D)$ associated with
a compact foliated Riemannian manifold $(M,{\cal F})$ as above,
the projection $P$ is equivariant with respect to the ${\Bbb R}$-actions 
on $S^*{\cal A}$ and 
$\bar{S}^{0}(G_{{\cal F}_N},{\cal L}(\pi^*E)\otimes\Omega^{1/2})$
given by $\alpha_t$ and $G^*_t$ respectively:
$$
P(\alpha_t(b))=G^*_t(P(b)), b\in S^*{\cal A}, t\in {\Bbb R}.
$$
\et
The organization of the paper is as follows.

Section~\ref{wave:sect} is devoted to investigation of the transversal wave
equation: norm estimates and a description in terms of Fourier integral 
operators. In particular, it contains a proof of Theorem~\ref{wavelapl} 
(see Theorem~\ref{waveform}).

In Section~\ref{sing}, we prove Theorem~\ref{DG} on singularities 
of the Fourier transform
$\theta(t)$ of the spectrum distribution function of a transversally 
elliptic operator (see Theorem~\ref{Poisson}).

Section~\ref{Eg:sect} is devoted to a proof of Theorem~\ref{Egorov0}
(see Theorem~\ref{Egorov}).

In Section~\ref{ncg}, we discuss the noncommutative geodesic flow
in transversal geometry of foliations and prove Theorem~\ref{noncom:geo} 
(see Theorem~\ref{noncom:flow}).
 
Section~\ref{short} is devoted to the short exact sequence
in transversal pseudodifferential calculus on a foliated manifold.
This is a technical result, which we need 
in Section~\ref{ncg}, but we believe that it is of independent interest. 

Finally, in Section~\ref{geom}, we discuss in more details geometrical 
examples of transversally elliptic
operators in Riemannian foliations and related objects.

\begin{ack} The work was partially done during visits to the
Humboldt University in Berlin supported by the Volkswagen Stiftung 
and to the Forschungsinstitut f\"ur Mathematik, ETH Z\"urich, and I wish 
to express my gratutide to them for hospitality and support. I would like
to thank E.Leichtnam for helpful discussions and his interest to this work.
\end{ack}

\section{Wave group}
\label{wave:sect}
Throughout in the paper, $(M,{\cal F})$ is a compact foliated manifold,
and $E$ is a Hermitian vector bundle on $M$. We will use the following
notation: $N^*{\cal F}$ is the conormal bundle to ${\cal F}$,
$G$ is the holonomy groupoid of ${\cal F}$. 
Let us briefly recall its definition. Let $\sim_h$ be an
equivalence relation on the set of leafwise paths $\gamma:[0,1]
\rightarrow M$, setting $\gamma_1\sim_h \gamma_2$ if $\gamma_1$ and 
$\gamma_2$ have the same initial 
and final points and the same holonomy maps. The holonomy groupoid 
$G$ is the set of $\sim_h$ equivalence classes of leafwise
paths. $G$ is equipped with the source and the range
maps $s,r:G\rightarrow M$ defined by $s(\gamma)=\gamma(0)$ and 
$r(\gamma)=\gamma(1)$. We will identify a point $x\in  M$ with
an element of $G$ given by the corresponding constant path: 
$\gamma(t)=x, t\in [0,1]$. Recall also that, for any $x\in M$, the set 
$G_{\cal F}^x=\{\gamma\in G : r(\gamma)=x\}$ is the covering
of the leaf through the point $x$, associated with the holonomy group 
of the leaf.  

We assume that the vector bundle $E$ is holonomy equivariant,
that is, it is equipped with an isometrical action
$$
T(\gamma):E_x\rightarrow E_y,\ \gamma\in G, \gamma:x\rightarrow y
$$
of $G$ in fibres of $E$.
We will use the classes $\Psi^{m,-\infty}(M,{\cal F},E)$ of transversal
pseudodifferential operators and the anisotropic Sobolev spaces
$H^{s,t}(M,{\cal F},E)$ (see \cite{noncom} for the definitions).

\begin{defn} 
\label{op}
Let us say that an operator $R$,
acting in $C^{\infty}(M,E)$,
belongs to the class ${\rm OP}^{l,-\infty}(M,N^*{\cal F},E)$, if, for any
$s,t,\alpha\in {\Bbb R}$, $R$ defines a continuous map
$$
R:H^{s,t}(M,{\cal F},E)\rightarrow H^{s+\alpha, t-\alpha-l}(M,{\cal F},E).
$$
\end{defn}

\begin{rem} 
\label{infty}
Any operator $R\in
\Psi^l(M,E)$, which has order $-\infty$ in some neighborhood of
$N^*{\cal F}$, belongs to the class ${\rm OP}^{l,-\infty}(M,N^*{\cal F},E)$
\cite{noncom}.
\end{rem}
\begin{rem}
\label{smooth}
Using mapping properties of operators from $\Psi^{m,-\infty}(M,{\cal F},E)$
 in $H^{s,t}(M,{\cal F},E)$ \cite{noncom}, it is easy to see that, 
if $A\in {\rm OP}^{l,-\infty}(M,N^*{\cal F},E)$ and 
$B\in \Psi^{m,-\infty}(M,{\cal F},E)$, then $AB$ and $BA$ are in 
$\Psi^{-\infty}(M,E)$.    
\end{rem}

For a pseudodifferential operator $P$, acting in $C^{\infty}(M,E)$,
let $p$ denote the principal symbol of $P$, and $\sigma_P$
the transversal principal symbol of $P$, which is the restriction
of $p$ to the punctured conormal bundle
$\tilde{N}^*{\cal F}=N^*{\cal F}\setminus\{0\}$.

Throughout in this paper, we will assume that
$P$ is an operator, acting in
$C^{\infty}(M,E)$, which satisfies the following conditions:

(H1) $P$ has the form
$$P=P_1+R_1,$$
where:
\medskip\par
(a) $P_1\in \Psi^1(M,E)$ is a transversally elliptic operator with the
scalar principal symbol and the holonomy invariant, positive 
transversal principal symbol;

(b) $R_1\in {\rm OP}^{1,-\infty}(M,N^*{\cal F},E)$.
\medskip\par
(H2) $P$ is an essentially self-adjoint operator in $L^2(M,E)$ (with the 
domain $C^{\infty}(M,E)$).

\begin{rem}
If $A\in \Psi^m(M,E)$ is a transversally elliptic
operator with the scalar principal
symbol and the holonomy invariant, positive transversal principal 
symbol, which is a self-adjoint positive operator in $L^2(M,E)$, then, by 
\cite{noncom} the operator
$P=A^{1/m}$ satisfies the conditions (H1) and (H2). 

Applying this statement to the transversal Laplacian, our basic geometric
example, we get that the operator $|D_H|$ satisfies the conditions (H1) 
and (H2). We refer the reader to Section~\ref{geom} for a more detailed 
discussion of notions connected with the operator $|D_H|$.  
\end{rem}

\bl
Any operator $P$, satisfying the conditions (H1) and (H2), can be represented
in the form
\be
\label{elliptic}
P=P_2+R_2,
\ee
where:
\medskip\par
(a) $P_2\in \Psi^1(M,E)$ is an essentially self-adjoint,
elliptic operator with the
positive scalar principal symbol and the holonomy invariant transversal
principal symbol such that the complete symbols of $P_1$ and $P_2$
coincides up to infinite order in some neighborhood of $N^*{\cal F}$;
\medskip\par
(b) $R_2\in {\rm OP}^{1,-\infty}(M,N^*{\cal F},E)$.
\el

\p Let $p_1\in S^1(I^n\times {\Bbb R}^n, {\cal L}
({\Bbb C}^N))$ be the
complete symbol of the operator $P_1$ in some foliated coordinate chart.
Assume that $p_1(x,y,\xi,\eta)$ is invertible for any
$(x,y,\xi,\eta)\in U, |\xi|^2+|\eta|^2>R^2$, where $R>0$, $U$ is a
conical neighborhood of $\eta=0$.
Take any function $\phi\in C^{\infty}(I^n\times {\Bbb R}^n),
\phi=\phi(x,y,\xi,\eta), x\in I^p, y\in I^q, \xi \in {\Bbb R}^p,
\eta \in {\Bbb R}^q$, homogeneous of degree $0$ in $(\xi,\eta)$
for $|\xi|^2+|\eta|^2>1$, which is supported in some conical neighborhood
of $\eta=0$ and is equal to $1$ in $U$, and put
$$
p_2(x,y,\xi,\eta)=\phi p_1(x,y,\xi,\eta) + (1-\phi)
(1+|\xi|^2+|\eta|^2)^{1/2}.
$$
Take $P_2$ to be the operator $p_2(x,y,D_x,D_y)$
with the complete symbol $p_2$ (or, more precisely,
$p_2(x,y,D_x,D_y)+p_2(x,y,D_x,D_y)^*$ to provide self-adjointness)
and $R_2=P-P_2$. By Remark~\ref{infty}, $P_1-P_2\in 
{\rm OP}^{1,-\infty}(M,N^*{\cal F},E)$, that completes immediately the proof.
\medskip
\par
By the spectral theorem, any self-adjoint operator $P$ in $L^2(M,E)$ defines
a strongly continuous semigroup $e^{itP}$ of bounded operators in $L^2(M,E)$.
Our next goal is to state mapping properties of operators $e^{itP}$ for
an operator $P$ under the conditions (H1) and (H2).

\begin{defn} 
Let us say that an operator
$K: C^{\infty}(M,E)\rightarrow C^{\infty}(M,E)$ belongs to the class
${\rm OP}^{m,\mu}(M,{\cal F},E)$, if, for any $s,k\in {\Bbb R}$, $K$ defines
a continuous map
$$
K: H^{s,k}(M,{\cal F},E)\rightarrow H^{s-m,k-\mu}(M,{\cal F},E).
$$
\end{defn}

\bp
\label{mapping}
Given an operator $P$, satisfying (H1) and (H2), the operator
$e^{itP}$ belongs to the class ${\rm OP}^{0,0}(M,{\cal F},E)+
{\rm OP}^{2,-\infty}(M,N^*{\cal F},E)$.
\ep

\p
Represent $P$ as in (\ref{elliptic}):
$P=P_2+R_2$.
By the usual Duhamel formula, we have
\be
\label{Duhamel}
e^{itP}=W(t)+R(t),
\ee
where $W(t)=e^{itP_2}$ and 
\be
\label{RT}
R(t)=\int_0^t e^{i\tau P_2}\,R_2\,e^{i(t-\tau)P}\,d\tau.
\ee

Let us check that $W(t)\in {\rm OP}^{0,0}(M,{\cal F},E)$.
Since $P_2$ is elliptic, $W(t)$ is, clearly, bounded
as an operator from $H^{s,0}(M,{\cal F},E)$ to $H^{s,0}(M,{\cal F},E)$ 
for any $s$. Let $\Lambda_{0,k}$ be a tangentially elliptic pseudodifferential
operator with the positive principal symbol, and let $\Lambda_{0,-k}$
be a leafwise parametrix of $\Lambda_{0,k}$. Then,
by holonomy invariance of the transversal principal symbol of $P_2$,
$$
\Lambda_{0,k}P_2\Lambda_{0,-k}=P_2+K_k,
$$
where $K_k\in \Psi^{0,0}(M,{\cal F},E)$, and, therefore, is bounded
in $H^s(M,E)$ for any $s$. This yields boundedness of $W(t)$
in the general case of the space $H^{s,k}(M,{\cal F},E)$ with $k\not= 0$.

By essential self-adjointness of $P$, the  operator
$e^{itP}$ is bounded in $L^2(M,E)$, and, using (\ref{RT}) and mapping 
properties of $e^{itP_2}$ and $R_2$, it is easy to see 
 $R(t)$ defines a bounded map from $L^2(M,E)$ to $H^{s,-s-1}(M,{\cal F},E)$ 
for any $s$, and, by duality, from $H^{s,-s+1}(M,{\cal F},E)$ to 
$L^2(M,E)$. From this, we get easily that 
$R(t)\in {\rm OP}^{2,-\infty}(M,N^*{\cal F},E)$,
if we write (\ref{RT}) as 
$$
R(t)=\int_0^t e^{i\tau P_2}\,R_2\,(e^{i(t-\tau)P_2}+R(t-\tau))\,d\tau.
$$

Using the standard Fourier integral representation for the operator
$W(t)$ of Proposition~\ref{mapping} (see, for instance, \cite{H4}),
we immediately obtain the following
description of the operator $e^{itP}$ in terms of Fourier integral
operators (it yields Theorem~\ref{wavelapl}, if applied to $|D_H|$).

\bt
\label{waveform}
Let $P$ be an operator, satisfying the conditions (H1) and (H2).
Let ${\tilde p}\in S^1(T^*M)$ be any scalar elliptic symbol, which 
is equal to the principal symbol $p$ of $P_1$ in some conical neighborhood
of $N^*{\cal F}$.
Let $\tilde{f}_t$ be the Hamiltonian flow of $\tilde{p}$.
For any $t\in {\Bbb R}$, let $\Lambda_{\tilde p}(t)$ be
the graph of the symplectomorphism $\tilde{f}_t$
in $\tilde{T}^*M=T^*M\setminus\{0\}$:
$$
\Lambda_{\tilde p}(t)=\{((x,\xi),(y,\eta))\in 
\tilde{T}^*M\times \tilde{T}^*M : (x,\xi)=\tilde{f}_t(y,\eta)\}.
$$
Then the operator $e^{itP}$ has the form
\begin{equation}
\label{wave}
e^{itP}=W(t)+R(t),
\end{equation}
where $W(t)$ is a Fourier integral operator associated with the canonical
relation $\Lambda_{\tilde p}(t)$, $W(t)\in I^0(M\times M, 
\Lambda'_{\tilde{p}}(t);\Omega^{1/2}(M\times M))$,
$R(t)$ is a smooth family of operators from 
${\rm OP}^{2,-\infty}(M,N^*{\cal F},E)$.
\et

\section{Singularities of the Fourier transform of the spectral function}
\label{sing}
In this Section, using Theorem~\ref{waveform} and the composition theorem 
for Fourier integral operators, we get
a result on the singularities of the Fourier transform
$\theta(t)$ of the
distributional spectrum distribution function of an operator $P$ under
the conditions (H1) and (H2).

Recall that, for any holonomy equivariant vector bundle $E$, we have
a $\ast$-representation $R_E$ of the involutive
algebra $C^{\infty}_c(G)$ in $L^2(M,E)$. For any
$k\in C^{\infty}_c(G)$, the operator $R_E(k)$
is given by the formula
\be
\label{repr}
R_E(k)u(x) =\int_{ G^{x}}  k(\gamma ) T(\gamma)[u(s(\gamma ))] 
d\lambda ^{x}(\gamma ), x\in  M, u\in C^{\infty}(M,E),
\ee
where $\lambda^x$ is the Riemannian volume on $G^x$.
 
A function $\theta(t)\in {\cal D}'(G), t\in {\Bbb R},$ is defined 
by the formula
$$
\theta(t)= {\rm tr}_G\; e^{itP}, 
$$
where ${\rm tr}_G$ is the $G$-trace functional introduced in
\cite{noncom}, or, equivalently,
$$
<\theta(t),k>= \hbox{tr}\; R_E(k)e^{itP}, k\in C^{\infty}_c(G).
$$
If $N(\lambda)\in {\cal D}'(G)$ denotes the spectrum distribution
function of the operator $P$, $<N(\lambda),k>={\rm tr}\; R_E(k)E(\lambda),
\lambda\in {\Bbb R}$, then $\theta(t)$ is the Fourier transform of 
$N(\lambda)$. 

\bt
\label{Poisson}
Given an operator $P$ under the conditions (H1) and (H2), for any $k\in
C^{\infty}_c(G)$, the function $\theta_k(t)=<\theta(t),k>$ is smooth
outside the set of $t\in {\Bbb R}$ such that there exists $\gamma
\in \supp k, \nu\in \tilde{N}^*_{r(\gamma)}{\cal F}$, with
$$
\tilde{f}_{-t}\ dh_{\gamma}^*(\nu)=\nu,
$$
where $dh_{\gamma}$ is the linear holonomy map of $\gamma$.
\et

\p
Let $k\in C^{\infty}_c(G)$. By Theorem~\ref{waveform}, we have
$$ 
R_E(k)e^{itP}=R_E(k)W(t)+R_E(k)R(t).
$$
By Remark~\ref{smooth}, $R_E(k)R(t)\in \Psi^{-\infty}(M,E)$. To study
the operator $R_E(k)W(t)$, we use the following description of
$R_E(k)$ as a Fourier integral operator \cite{noncom}.
  
Let ${\cal F}_N$ be the foliation in $\tilde{N}^*{\cal F}$,
which is the lift of ${\cal F}$
via a leafwise flat (Bott) connection in $\tilde{N}^*{\cal F}$ \cite{Molino},
and $G_{{\cal F}_N}$ the holonomy groupoid  of ${\cal F}_N$.
The holonomy groupoid $G_{{\cal F}_N}$ consists of all triples
$(\gamma,\nu)\in G_{\cal F}\times H^{*}$ such
that $r(\gamma)=\pi(\nu)$, and $s(\gamma)=dh_{\gamma}^{*}(\nu)$ 
with the source map $s:G_{{\cal F}_H}\rightarrow
H^{*}, s(\gamma,\nu)=dh_{\gamma}^{*}(\nu)$ and the range map
$r:G_{{\cal F}_H}\rightarrow
H^{*}, r(\gamma,\nu)=\nu$. We have a map
$\pi_G:G_{{\cal F}_H}\rightarrow G_{\cal F}$, 
given by $\pi_G(\gamma,\nu)=\gamma$.

For any $k\in C^{\infty}_c(G)$, the operator $R_E(k)$ belongs to
$\Psi^{0,-\infty}(M,{\cal F},E)$.
Operators from $\Psi^{0,-\infty}(M,{\cal F},E)$ can be considered
as Fourier integral operators, associated with an immersed
canonical relation, which is the image of $G_{{\cal F}_N}$ under
the mapping
$$
(r,s):G_{{\cal F}_N}\rightarrow T^*M\times T^*M,
$$
given by the source and the target mappings of the groupoid $G_{{\cal F}_N}$, 
$R_E(k)\in I^{-p/2}(M\times M,G'_{{\cal F}_N};\Omega^{1/2}(M\times M))$.

Now we apply the composition theorem for Fourier integral operators.
As in Theorem~\ref{waveform}, let ${\tilde p}\in S^1(T^*M)$ be 
any scalar elliptic symbol, which is equal to
the principal symbol $p$ of $P_1$ in some conical neighborhood
of $N^*{\cal F}$. By holonomy invariance of the transversal
principal symbol of $P$, the Hamiltonian flow $\tilde{f}_t$ of $\tilde{p}$
leaves $N^*{\cal F}$ invariant, and we denote by 
$f_t$ the restriction of the flow $\tilde{f}_t$ to $N^*{\cal F}$.
It is easy to see that $f_t$ doesn't depend on a choice of $\tilde{p}$. 

Since $\Lambda'_{\tilde{p}}(t)$ and $G'_{{\cal F}_N}$ intersects 
transversally, by the composition theorem for Fourier integral operators
\cite{H4},
we get that the operator $R_E(k)e^{itP}$ is a Fourier integral operator,
associated with the immersed canonical relation $\Lambda'_{\tilde{p}}(t)
\circ G'_{{\cal F}_N}$, $R_E(k)e^{itP}\in I^{-p/2}(M\times M, 
\Lambda'_{\tilde{p}}(t)\circ G'_{{\cal F}_N};\Omega^{1/2}(M\times M))$.
It is easy to see that $\Lambda_{\tilde{p}}(t)\circ G_{{\cal F}_N}$ 
is the image of $G_{{\cal F}_N}$ under the mapping
$$
(r, f_{-t}\circ s):G_{{\cal F}_N}\rightarrow T^*M\times T^*M,
$$
where $f_{-t}\circ s$ is the composition of the source map
$s :G_{{\cal F}_N}\rightarrow N^*{\cal F}$ and the diffeomorphism
$f_{-t}$ of $N^*{\cal F}$ to itself.

Now the proof of the theorem is completed by the same arguments as in 
\cite{DG}. Using the functorial properties of the wave-fronts sets, 
we state that the singularities of $\theta_k(t)$ are 
contained in the set of all $t$ such that the intersection of  
$\Lambda_{\tilde{p}}(t)\circ G_{{\cal F}_N}$ with the diagonal
in $T^*M\times T^*M$ is not empty. Due to above description of 
$\Lambda_{\tilde{p}}(t)\circ G_{{\cal F}_N}$, it is equivalent to the 
statement of the theorem.

\begin{rem}
Just as in \cite{DG} for the case of an elliptic operator on a compact
manifold, we can prove that the distribution $\theta_k$ is Lagrangian and
get a Duistermaat-Guillemin type trace formula for transversally elliptic
operators. We hope to discuss this in a forthcoming paper.
\end{rem}

\section{Egorov type theorem}
\label{Eg:sect}
As above, let $P$ be an operator, satisfying the conditions (H1) and (H2).
For any operator $K: C^{\infty}(M,E)\rightarrow C^{\infty}(M,E)$,
we define an operator $\Phi_t^P(K): C^{\infty}(M,E)\rightarrow
C^{\infty}(M,E)$ by the formula
\begin{equation}
\label{phi}
\Phi_t^P(K)=e^{itP}Ke^{-itP}.
\end{equation}
It is clear that the formula (\ref{phi}) defines a map
$$
\Phi_t^P:{\cal B}(L^2(M,E))\rightarrow {\cal B}(L^2(M,E)).
$$
Proposition \ref{mapping} implies immediately mapping properties
of operators $\Phi_t^P(K)$ in Sobolev spaces $H^{s,k}(M,{\cal F},E)$.

\bp
Given an operator $P$ satisfying the conditions (H1) and (H2) and
$K\in {\rm OP}^{m,\mu}(M,{\cal F},E)$, we have
$$
\Phi_t^P(K)\in {\rm OP}^{m,\mu}(M,{\cal F},E)
+{\rm OP}^{m+\mu+4,-\infty}(M,N^*{\cal F},E).
$$
\ep

Now we turn to a description of the operator $\Phi_t^P(K)$ when
$K$ is a transversal pseudodifferential operator.
Let $\Omega^{1/2}({\cal F}_N)$ be
the vector bundle of leafwise half-densities on $N^*{\cal F}$. 
We can lift this bundle to vector bundles
$s^*(\Omega^{1/2}({\cal F}_N))$ and $r^*(\Omega^{1/2}({\cal F}_N))$ 
on $G_{{\cal F}_N}$ via the source and target mappings $s$ and 
$r$ respectively and form a bundle $\Omega^{1/2}$ on	$G_{{\cal F}_N}$
by
$$
\Omega^{1/2}=r^*(\Omega^{1/2}({\cal F}_N))\otimes
(s^*(\Omega^{1/2}({\cal F}_N)))^*.
$$

The space $S^{m}(G_{{\cal F}_N},{\cal L}(\pi^*E)\otimes\Omega^{1/2})$
of transversal symbols is the space of all
$s\in C^{\infty}(G_{{\cal F}_N},{\cal L}(\pi^*E)\otimes\Omega^{1/2})$
homogeneous of degree $m$ such that $\pi_G({\rm supp}\; s)$ is compact
in $G_{\cal F}$.
The space
$$
S^{*}(G_{{\cal F}_N}, {\cal L}(\pi^*E)\otimes \Omega^{1/2})=
\bigcup_mS^m(G_{{\cal F}_N}, {\cal L}(\pi^*E)\otimes \Omega^{1/2})
$$
carries the structure of an involutive algebra, defined by
its embedding into the foliation algebra
$C^{\infty}_c(G_{{\cal F}_N}, {\cal L}(\pi^*E))$.

By \cite{noncom}, there is the half-density principal symbol mapping
$$
\Psi^{m,-\infty}(M,{\cal F},E)\rightarrow
S^{m}(G_{{\cal F}_N},{\cal L}(\pi^*E)\otimes\Omega^{1/2}),
$$
which is an $\ast$-isomorphism.

Let ${\tilde p}\in S^1(T^*M)$ be any scalar elliptic symbol, which 
is equal to the principal symbol $p$ of $P_1$ in some conical neighborhood
of $N^*{\cal F}$. Let $\tilde{f}_t$ be the Hamiltonian flow of $\tilde{p}$,
and $\{f_t\}$ its restriction to  $N^*{\cal F}$. 
Holonomy invariance of the transversal principal symbol $\sigma_P$ 
implies existence of a flow $\{F_t\}$ on $G_{{\cal F}_N}$ such that
$$ F_t\circ s=F_t\circ r=f_t.$$

\begin{defn} 
Given an operator $P$ under the conditions (H1) and (H2), 
{\bf the transversal bicharacteristic flow} of $P$ is
the one-parameter group $F^*_t$ of automorphisms of the
involutive algebra 
$S^m(G_{{\cal F}_N},{\cal L}(\pi^*E)\otimes\Omega^{1/2})$,
induced by $F_t$.
\end{defn}

\begin{rem}
It can be easily seen from Theorem~\ref{Egorov} that the definition of 
the transversal bicharacteristic flow  
doesn't depend on a choice of a representation of the operator $P$ 
in the form (H1).  
\end{rem} 

\noindent The following is Theorem~\ref{Egorov0} of Introduction.
\bt
\label{Egorov}
Given an operator $K\in \Psi^{m,-\infty}(M,{\cal F},E)$ with the principal
symbol $k\in S^{m}(G_{{\cal F}_N},{\cal L}(\pi^*E)\otimes\Omega^{1/2})$,
the operator $\Phi_t^P(K)$ is an operator of class
$\Psi^{m,-\infty}(M,{\cal F},E)$ with the principal symbol
$k(t)\in S^m(G_{{\cal F}_N},{\cal L}(\pi^*E)\otimes\Omega^{1/2})$,
given by
$$
k(t)=F^*_t(k).
$$
\et

\p By (\ref{wave}), we have
$$
\Phi_t^P(K)=P(t)+C(t),
$$
where
$$
P(t)=W(t)KW(-t), C(t)=R(t)KW(-t)+W(t)KR(-t)+R(t)KR(-t).
$$
By Proposition~\ref{mapping}, the operator $C(t)$ belongs to
$\Psi^{-\infty}(M,E)$. The description of $P(t)$ as an operator
from $\Psi^{m,-\infty}(M,{\cal F},E)$ and the formula for its
principal symbol follows straightforwardly from the composition 
theorem of Fourier integral operators (see, for instance, \cite{H4}).

\begin{rem} 
We can also write down a differential equation on
the principal symbol
$k(t)\in S^m(G_{{\cal F}_N},{\cal L}(\pi^*E)\otimes\Omega^{1/2})$
as follows.
Let $H_p$ be a vector field on $N^*{\cal F}$, which is the restriction
of the Hamiltonian vector field of $p$. Define a vector field $v_p$
on $G_{{\cal F}_N}$ as
$$
v_p(\gamma,\eta)=dr^*(H_p(y,\eta))+ds^*(H_p(x,\eta)),
\gamma:(x,\eta)\rightarrow (y,\eta),
$$
where $dr^*(H_p(y,\eta))$, $ds^*(H_p(x,\eta))$ are the lifts of the
vector field $H_p$ via the covering maps
\ba
s: G^{(y,\eta)}_{{\cal F}_N}&=&\{\gamma'\in G_{{\cal F}_N}: 
r(\gamma',\eta)=(y,\eta)\}\rightarrow N^*{\cal F},\\
r: (G_{{\cal F}_N})_{(x,\eta)}&=&\{\gamma'\in G_{{\cal F}_N}: 
s(\gamma',\eta)=(x,\eta)\}
\rightarrow N^*{\cal F}.
\ea
Then $k(t)$ is a solution of the following Cauchy problem
\be
\label{diff}
\frac{dk(t)}{dt}={\cal L}_{v_p}(k(t)), t>0, k(0)=k,
\ee
where ${\cal L}_{v_p}$ denotes the Lie derivative.

It would be very interesting to give an interpretation of the 
equation~(\ref{diff}) as a noncommutative Hamiltonian equation. 
 
\end{rem}

\section{Noncommutative geodesic flow in transversal geometry}
\label{ncg}
In this Section, we will give an interpretation of the Egorov
type theorem, Theorem~\ref{Egorov}, in terms of the
corresponding noncommutative geodesic flow.
As in \cite{noncom}, let us consider spectral triples
$({\cal A},{\cal H},D)$ associated with a compact foliated Riemannian
manifold $(M,{\cal F})$ of the form:
\begin{enumerate}
\item The involutive
algebra ${\cal A}$ is the algebra $\cinf_c(G)$;
\item The Hilbert space ${\cal H}$
is the space $L^2(M,E)$ of $L^2$-sections of a holonomy
equivariant Hermitian vector bundle $E$, on which
an element $k$ of the algebra ${\cal A}$ is represented via
the $\ast$-representation $R_E$;
\item The operator $D$ is a first order self-adjoint
transversally elliptic operator with the holonomy invariant transversal
principal symbol such that the operator $D^2$ has the scalar principal symbol
and self-adjoint.
\end{enumerate}

Let $\delta$ be a derivative on ${\cal L}({\cal H})$ given by 
$\delta(T)=[|D|,T], T\in {\cal L}({\cal H})$. By \cite{noncom}, ${\cal A}$
satisfies the following smoothness condition: for any $a\in {\cal A}$, 
$a$ and $[D,a]$ belong to the domain of $\delta^n$ for any $n$. 
Let ${\cal B}$ be the algebra generated by $\delta^n(a), a\in {\cal A}, 
n\in {\Bbb N}$. As shown in \cite {noncom}, ${\cal B}$ is contained in
$\Psi^{0,-\infty}(M,{\cal F},E)$. 
Recall \cite{Co-M,Sp-view} that ${\rm OP}^{\alpha}$ denotes 
the space of operators in ${\cal H}$ of order ${\alpha}$, that means that
$P\in {\rm OP}^{\alpha}$ iff $P|D|^{-\alpha}\in \bigcap_n {\rm Dom}\;\delta^n$.
Unlike the case  ${\cal A}$ is an unital algebra as in \cite{Sp-view},
the corresponding definition of the algebra $\Psi^*({\cal A})$ 
of pseudodifferential operators is not convenient in our case,  
because it may happen that operators from $\Psi^*({\cal A})$ are not
smoothing in leafwise directions. Therefore, we should
modify the definition of the algebra $\Psi^*({\cal A})$ of \cite{Co-M,Sp-view} 
if we wish to get a natural generalization of the algebra of 
pseudodifferential operators in a nonunital case. Since we are only interested
in the corresponding bicharacteristic flow here, it is enough for us to 
give an appropriate definition of the norm closure of 
$\Psi^*({\cal A})$ that can be done rather roughly.
  
Define the algebra $\Psi^*_f({\cal A})$ of pseudodifferential operators 
as the set of operators, represented as a finite sum:
\be
\label{psif}
P= b_q|D|^q + b_{q-1}|D|^{q-1}+\ldots+b_{-N}|D|^{-N}, b_j\in {\cal B}.
\ee 
Now let ${\cal C}_f$ be the algebra 
$${\cal C}={\rm OP}^0\bigcap \Psi^*_f({\cal A}),$$
and $\bar{\cal C}$ be the closure of ${\cal C}$ in ${\cal L}({\cal H})$.

Finally, we come to the definition of the unitary cotangent bundle.
\begin{defn}
Given a spectral triple $({\cal A}, {\cal H}, D)$,
{\bf the unitary cotangent bundle} $S^*{\cal A}$ is defined
as the quotient of the $C^*$-algebra ${\bar{\cal C}}$
by the ideal ${\cal K}$ of compact operators in ${\cal H}$.
\end{defn} 

Now we pass to a description of the objects introduced above for the
spectral triples associated with a Riemannian foliation.
We start with pseudodifferential calculus.

\bl
\label{Psi}
Let $({\cal A},{\cal H},D)$ be a spectral triple associated with
a compact foliated Riemannian manifold $(M,{\cal F})$ as above.
The algebra $\Psi^*_f({\cal A})$ is contained in the 
$C^*$-algebra $\Psi^{0,-\infty}(M,{\cal F},E)$.
\el
\p
Let $P\in \Psi^*_f({\cal A})$ as in (\ref{psif}).  
By \cite{noncom}, we know that the operator $|D|^j, j\in {\Bbb Z},$ 
has the form
\begin{equation}
\label{decom}
|D|^j= P(j)+R(j),
\end{equation}
where $P(j)\in \Psi^{j}(M,E)$ has the principal symbol $p_{j}$, 
being equal to $|a_1|^j$ in some conic neighborhood of $N^*{\cal F}$
($a_1$ is the principal symbol of $D$),
and, for any $s,l,t,j$, $R(j)$ defines a continuous mapping
\begin{eqnarray}
\label{Tj}
R(j):H^{s,l}(M,{\cal F},E)&\rightarrow &H^{t,s+l-t-j}(M,{\cal F},E), j>0,
\nonumber\\
R(j):H^{s,l}(M,{\cal F},E)&\rightarrow &H^{t,s+l-t}(M,{\cal F},E), j\leq 0.
\end{eqnarray}

By (\ref{decom}), we have $b_j|D|^j=b_jP(j)+b_jR(j)$, where 
$b_jP(j)\in \Psi^{j,-\infty}(M,{\cal F},E)$ has the principal symbol, 
being equal to $b_j|a_1|^j$, and, by (\ref{Tj}) and mapping properties
of operators from $\Psi^{0,-\infty}(M,{\cal F},E)$ \cite{noncom}, the
operator $b_jR(j)$ belongs to $\Psi^{-\infty}(M,E)$. So we conclude that
the operator $b_j|D|^j$ belongs to $\Psi^{j,-\infty}(M,{\cal F},E)$ 
for any $j$, that completes the proof.
\medskip
\par
To provide a description of the unitary cotangent bundle, we need the
short exact sequence in transversal pseudodifferential calculus
on a Riemannian foliation. Here we will formulate the corresponding
results, referring the reader to Section~\ref{short} for proofs.

Let $\bar{\Psi}^{0,-\infty}(M,{\cal F},E)$ be the closure  of
$\Psi^{0,-\infty}(M,{\cal F},E)$ in ${\cal L}(L^2(M,E))$ and
$\bar{S}^{0}(G_{{\cal F}_N},{\cal L}(\pi^*E)\otimes\Omega^{1/2})$
the closure  of
$S^{0}(G_{{\cal F}_N},{\cal L}(\pi^*E)\otimes\Omega^{1/2})$
in the uniform operator topology of $L^2(G_{{\F}_N},s^*(\pi^*E))$. 
By Proposition~\ref{cont}, the symbol map $\sigma$ extends by continuity 
to a map 
$$
\bar{\sigma}:\bar{\Psi}^{0,-\infty}(M,{\cal F},E)
\rightarrow
\bar{S}^{0}(G_{{\cal F}_N},{\cal L}(\pi^*E)\otimes\Omega^{1/2}).
$$
Denote by $I_{\sigma}$ the kernel of the symbol mapping $\bar{\sigma}$.
It is a two-sided ideal in the $C^*$-algebra
$\bar{\Psi}^{0,-\infty}(M,{\cal F},E)$. Let us identify elements
of $\bar{\Psi}^{0,-\infty}(M,{\cal F},E)$ with corresponding bounded
operators in $L^2(M,E)$. 

\bp
\label{seq}
(1) There is defined a short exact sequence
$$
0\rightarrow I_{\sigma}\rightarrow \bar{\Psi}^{0,-\infty}(M,{\cal F},E)
\rightarrow
\bar{S}^{0}(G_{{\cal F}_N},{\cal L}(\pi^*E)\otimes\Omega^{1/2}).
$$

(2) The ideal $I_{\sigma}$ contains the ideal ${\cal K}$
of compact operators in $L^2(M,E)$:
\be
\label{incl}
{\cal K}\subset I_{\sigma}
\ee
\ep

By (\ref{incl}) and Lemma~\ref{Psi}, it is easy to see that 
the map $\bar{\sigma}$ defines a  surjective map
$$
P: S^*{\cal A}\rightarrow
\bar{S}^{0}(G_{{\cal F}_N},{\cal L}(\pi^*E)\otimes\Omega^{1/2}).
$$
\noindent The Egorov's type theorem, 
Theorem~\ref{Egorov}, provides a description
of the image of the noncommutative geodesic flow in $S^*{\cal A}$,
defined by the triple $({\cal A}, {\cal H}, D)$, under this projection. 
First of all, recall the definition \cite{Sp-view}.

\begin{defn}
Given a spectral triple $({\cal A}, {\cal H}, D)$, {\bf the noncommutative
geodesic flow} is a one-parameter group $\alpha_t$ of automorphisms
of the algebra $S^*{\cal A}$, defined by
$$
\alpha_t(b)=e^{it|D|}be^{-it|D|}, b\in S^*{\cal A}, t\in {\Bbb R}.
$$
\end{defn}

Theorem~\ref{Egorov} can be reformulated in the following way 
(Theorem~\ref{noncom:geo} of Introduction).

\bt
\label{noncom:flow}
Let $({\cal A},{\cal H},D)$ be a spectral triple associated with
a compact foliated Riemannian manifold $(M,{\cal F})$ as above.
Then
$$
P(\alpha_t(b))=F^*_t(P(b)), b\in S^*{\cal A}, t\in {\Bbb R},
$$
where $F^*_t$ is the transversal bicharacteristic flow of the operator $|D|$,
which is a one-parameter group  of automorphisms of 
$\bar{S}^0(G_{{\cal F}_N},{\cal L}(\pi^*E)\otimes\Omega^{1/2})$.
\et

\section{The short exact sequence in transversal $\Psi$DO 
calculus} 
\label{short}
Here we state the short exact sequence
in transversal pseudodifferential calculus on a foliated manifold.
As above, let $(M,{\cal F})$ be a compact foliated manifold, and 
$E$ is an Hermitian vector bundle on $M$. It should 
be noted that we don't assume here the foliation ${\cal F}$ to be Riemannian,
and the bundle $E$ to be holonomy invariant. 
Let $\bar{\Psi}^{0,-\infty}(M,{\cal F},E)$ be the closure  of
$\Psi^{0,-\infty}(M,{\cal F},E)$ in the uniform operator topology
of $L^2(M,E)$, and
$\bar{S}^{0}(G_{{\cal F}_N},{\cal L}(\pi^*E)\otimes\Omega^{1/2})$
the closure  of
$S^{0}(G_{{\cal F}_N},{\cal L}(\pi^*E)\otimes\Omega^{1/2})$
in the uniform operator topology of 
$$
L^2(G_{{\F}_N},s^*(\pi^*E))=\bigcup_{\nu\in N^*\F}
L^2(G_{{\F}_N}^{\nu},s^*(\pi^*E)).
$$

\bp
\label{cont}
The principal symbol mapping
$$
\sigma:\Psi^{0,-\infty}(M,{\cal F},E)\rightarrow
S^{0}(G_{{\cal F}_N},{\cal L}(\pi^*E)\otimes\Omega^{1/2}),
$$
extends by continuity to a map
$$
\bar{\sigma}:\bar{\Psi}^{0,-\infty}(M,{\cal F},E)\rightarrow
\bar{S}^{0}(G_{{\cal F}_N},{\cal L}(\pi^*E)\otimes\Omega^{1/2}).
$$
\ep
\p It suffices to prove that
\be
\label{sigma}
\|\sigma(P)\|\leq \|P\|, P\in \bar{\Psi}^{0,-\infty}(M,{\cal F},E).
\ee
Take $P\in\Psi^{0,-\infty}(M,{\cal F},E)$ and 
$u\in \cinf_c(G_{{\F}_N}^{\nu},s^*(\pi^*E))$
for some $\nu\in \tilde{N}^*{\F}$ with $\supp u=K\in G_{{\F}_N}^{\nu}
\cong G^{\pi(\nu)}$ ($\pi : N^*{\F}\rightarrow M$ is the bundle map).

We will use the following theorem, which is a direct generalization
of the fundamental structure theorem for a neighborhood of a compact
leaf with finite holonomy and was used in \cite{thesis,Hurder} (this 
formulation is due to \cite{Hurder}).

Let $(V,{\cal V})$ be a compact foliated manifold.
Given a set $Z\subset V$ and $\varepsilon>0$ let $N(Z,\varepsilon)$
be the open neighborhood of $V$ consisiting of points which lie
within $\varepsilon$ of $Z$.
\bt
\label{reeb}
Let $L$ be a leaf in a compact foliated manifold $(V,{\cal V})$.
Given a compact subset $K\subset \tilde{L}_h$ and $\varepsilon>0$,
there exists $x\in I^q$ and $\delta>0$ so that for the open ball
$B(x,\delta)\subset I^q$ there is a foliated immersion
$\Pi: K\times B(x,\delta)\rightarrow V$ with
$\left. \Pi\right|_{K\times\{x\}}$, being the restriction to $K$ of the
covering map $\pi: \tilde{L}_h\rightarrow L$, and
$$
\Pi(K\times B(x,\delta))\subset N(\pi(K),\varepsilon).
$$
\et
Applying Theorem~\ref{reeb} in the case $(V,{\cal V})=
(N^*{\F},G_{\F_N})$ and $K=\supp u$,
we get a foliated immersion $\Pi: K\times B(Y,\delta)\times {\Bbb R}^q
\rightarrow N^*{\F}$
as in Theorem~\ref{reeb}. For any $(y,\eta)\in
B(Y,\delta)\times {\Bbb R}^q$, put
$$
\phi_{y,\eta}(y')=e^{i<y',\eta>}(\|\eta\|+1)^{q/4}
\chi((\|\eta\|+1)^{1/2}(y-y')), y'\in I^q,
$$
where $\chi\in\cinf_c({\Bbb R}^q)$ such that $\|\chi\|_2=1$, and, for
any $y\in B(Y,\delta), \lambda\in {\Bbb R}, \lambda\geq 1$, the function
$$
y'\mapsto \chi(\lambda(y-y'))
$$
is supported in $I^q$. Let $v_{y,\eta}\in C^{\infty}_c (G_{\F_N}), 
(y,\eta)\in B(Y,\delta)\times {\Bbb R}^q$ 
be defined on $\Pi(K\times B(Y,\delta))$ as
$$
v_{y,\eta}(\Pi(x,y'))=u(x)\phi_{y,\eta}(y'), x\in K, 
y'\in B(Y,\delta),
$$
and vanish outside of $\Pi(K\times B(Y,\delta))$.

Let us use the crossed product representation of the operator
$P$ in the foliated neighborhood $\Pi(K\times B(Y,\delta))$\cite{noncom}.
In the representation of $L^2(K\times B(Y,\delta),{\Bbb C}^r)$ as
the $L^2$ space of $L^2(B(Y,\delta), {\Bbb C}^r))$-valued
functions on $K$,
$L^2(K\times B(Y,\delta),{\Bbb C}^r)= L^2(K, L^2(B(Y,\delta),
{\Bbb C}^r))$,
we can write:
$$
P\bar{u}(x)=\int  P(x,x')\bar{u}(x')\,dx', x\in K,
$$
where $\bar{u}\in C_c^{\infty}(K, L^2(B(Y,\delta), {\Bbb C}^r)),
\bar{u}(x)(y)=u(x,y), x\in K, y\in B(Y,\delta)$ and,
for any $x\in K$, $x'\in K$, the operator $P(x,x')$ is a
pseudodifferential operator of order $0$ on $B(Y,\delta)$:
$$
P(x,x')v(y)=(2\pi)^{-q} \int  e^{i(y-y')\eta}k(x,x',y,\eta) v(y')\,dy'\,d\eta,
v\in C^{\infty}_c(B(Y,\delta),{\Bbb C}^r).
$$

The principal symbol $\sigma_P$
of $P$ is related with the principal symbols $\sigma_{P(x,x')}$
of the operators $P(x,x')$
as follows:
$$
\sigma_P(x,x',y,\eta)=\sigma_{P(x,x')}(y,\eta).
$$
By standard pseudodifferential calculus, we have
\ba
<\sigma_{P(x,x')}(y,\eta)\bar{u}(x'),\bar{u}(x)>_{L^2(B(Y,\delta), 
{\Bbb C}^r)}\\=
<\bar{u}\phi_{y,\eta}, P(x,x')(\bar{u}\phi_{y,\eta})>_{L^2(B(Y,\delta), 
{\Bbb C}^r)},
\ea
that implies
\be
\label{symbol2}
<\sigma(P)(y,\eta)u,u>_{L^2(K\times B(Y,\delta),{\Bbb C}^r)}
=<u\phi_{y,\eta}, P(u\phi_{y,\eta})>_{L^2(K\times B(Y,\delta),{\Bbb C}^r)}.
\ee
from where the desired estimate (\ref{sigma}) follows immediately.
\medskip\par
Now we complete the proof of Proposition~\ref{seq}:
\medskip\par
{\it Proof of Proposition~\ref{seq}.}\ 
(1) It remains only to prove that the map $\bar{\sigma}$ is surjective,
that can be easily done by means of an explicit construction of
a pseudodifferential operator with the given symbol.

(2) Since $u\phi_{y,\eta}\rightarrow 0$ weakly when $\eta\rightarrow
\infty$, by (\ref{symbol2}), we get  that $\sigma(P)=0$ for any
$P\in {\cal K}$.

\begin{rem}
\label{extens}
Let $C^{*}_E(G)$ be the closure of $R_E(C^{\infty}_c(G))$ 
in the uniform operator topology
of ${\cal L}(L^2(M,E))$ and $C^{\ast}_{r}(G)$ the reduced foliation 
$C^*$-algebra (see, for instance, \cite{F-Sk}).
By \cite{F-Sk}, we have the surjective  projection
$$
\pi_E : C^{\ast}_{E}(G) \rightarrow  C^{\ast}_{r}(G).
$$
Since $R_E(k)\in \Psi^{0,-\infty}(M,{\cal F},E)$ for any 
$k\in C^{\infty}_c(G)$, $C^{\ast}_{E}(G)$ is contained in  
$\bar{\Psi}^{0,-\infty}(M,{\cal F},E)$.
The corresponding set of principal symbols can be identified with $C^*_r(G)$.
So the principal symbol map $\bar{\sigma}$ provides 
an extension of $\pi_E$ to $\bar{\Psi}^{0,-\infty}(M,{\cal F},E)$. 
\end{rem}

\begin{rem}
It is easy to see that, for any $P\in \bar{\Psi}^{0,-\infty}(M,{\cal F},E)$, 
we have the inclusion
$$
{\rm spec}(\bar{\sigma}(P))\subset {\rm spec}(P),
$$
where ${\rm spec}(P)$ is the spectrum of $P$ in $L^2(M,E)$ and 
${\rm spec}(\bar{\sigma}(P))$ is the leafwise spectrum of 
$\bar{\sigma}(P)$ in $L^2(G_{{\F}_N},s^*(\pi^*E))$. 
Applied to operators $R_E(k)$ (see Remark~\ref{extens}), this fact provides
a proof of an ``easy'' part of  the spectral coindence theorem for 
tangentially elliptic operators (see \cite{tang} and references therein). 
As shown there, the inverse inclusion is not 
always true and is related with amenability of the foliation in question. 

A related remark is the following.
If  $I_{\sigma}={\cal K}$, then
the essential spectra of both operators, $P$ and $\bar{\sigma}(P)$, 
are the same:
$$
{\rm spec}_{\rm ess}(\bar{\sigma}(P))={\rm spec}_{\rm ess}(P), 
P\in \bar{\Psi}^{0,-\infty}(M,{\cal F},E).
$$
The spectrum of $\bar{\sigma}(P)$ is essential, 
${\rm spec}_{\rm ess}(\bar{\sigma}(P))=
{\rm spec}(\bar{\sigma}(P))$ (see \cite{tang} for the case $P=R_E(k)$),
but there are examples of operators $P$ with ${\rm spec}_{\rm ess}\;(P)\not=
{\rm spec}(P)$ \cite{adiab}.
\end{rem}

\section{Transversal geodesic flow}
\label{geom}
In this Section, we discuss geometrical examples of transversally elliptic
operators in Riemannian foliations.
As usual, we assume that $(M,{\cal F})$ is
a Riemannian foliation with a bundle-like metric $g_{M}$. Recall 
this means that $M$ satisfies one of the following equivalent
conditions (see \cite{Re}):
\begin{enumerate}
\item $(M,{\cal F})$ locally has the structure of a Riemannian
submersion;
\item Let $F=T{\cal F}$ be the tangent bundle to ${\cal F}$,
$H$ be the orthogonal complement to $F$, and $g_H$ be the restriction
 of $g_M$ to $H$. Then, for any $X\in F$, we have
$$
\nabla_X^{\cal F}g_H=0,
$$
where $\nabla^{\cal F}$ is the Bott connection on $H$;
\item the horizontal distribution $H$ is totally geodesic.
\end{enumerate}

There are two natural examples of geometric operators on Riemannian
foliations:
\medskip\par
(1) The Laplacian $\Delta_M$ of the Riemannian metric $g_M$. This
is a second order elliptic differential operator with the principal
symbol
$$
\sigma(\Delta_M)(x,\xi)=g_M(\xi,\xi)I_x, (x,\xi)\in \tilde{T}^*M.
$$
(2) The transversal (horizontal) Laplacian $\Delta_H$. We recall its
definition. The decomposition
$F\oplus H=TM$
induces a bigrading on $\Lambda T^{*}M$:
$$\Lambda^k T^{*}M=\bigoplus_{i=0}^{k}\Lambda^{i,k-i}T^{*}M,$$
where
$$\Lambda^{i,j}T^{*}M=\Lambda^{i}F^{*}\otimes
\Lambda^{j}H^{*}.$$
In this bigrading,
the de Rham differential $d$ can be written as
$$
d=d_F+d_H+\theta,
$$
where $d_F$ and $d_H$ are first order differential
operators, which are the tangential de Rham differential
and the transversal de Rham differential accordingly,
and $\theta$ is a zeroth order differential operator.

The transversal Laplacian is a second order elliptic differential operator
in the space $C^{\infty}(M,\Lambda T^{*}M)$:
$$
\Delta_H=d_H \delta_H+ \delta_Hd_H.
$$

\noindent Its principal symbol $\sigma(\Delta_H)$ is given by
$$
\sigma(\Delta_H)(x,\xi)=g_H(\xi,\xi)I_x, (x,\xi)\in \tilde{T}^*M.
$$
This formula can be rewritten as follows.
Let $P$ be a projection
$$P:T^{*}M \rightarrow T^{*}{\cal F},$$
which is adjoint to the natural inclusion
$T{\cal F}\rightarrow TM$. Then
$$
\sigma(\Delta_H)(x,\xi)=(g_M(\xi,\xi) - g_M(P\xi,P\xi))I_x,
(x,\xi)\in \tilde{T}^*M.
$$

\noindent Using this, it can be easily checked that
\begin{eqnarray}
\label{equal}
\sigma(\Delta_M)(x,\xi) &= &\sigma(\Delta_H)(x,\xi), (x,\xi)\in
N^*{\cal F},\nonumber\\
d\sigma(\Delta_M)(x,\xi) &= &d\sigma(\Delta_H)(x,\xi), (x,\xi)\in
N^*{\cal F}.
\end{eqnarray}
Let $p_2$ be any smooth function, which coincides with
$\sqrt{\sigma(\Delta_H)}$ in some conical neighborhood of
$N^*{\cal F}$.
The equalities (\ref{equal}) imply that the restrictions of the
Hamiltonian flows of $p_1=\sqrt{\sigma(\Delta_M)}$
and $p_2$ to $N^*{\cal F}$
coincide and are given by the restriction of the geodesic flow
$g_t$ of the Riemannian metric $g_M$ to $N^*{\cal F}$.
Thus, we see that the transversal bicharacteristic flows of the
operators $|D_H|=\sqrt{\Delta_H}$ and $\sqrt{\Delta_M}$ are 
given by the transversal part of the geodesic flow on $M$. 

Finally, if ${\cal F}$ is given by the fibres of a
Riemannian submersion $f:M\rightarrow B$, then there is a natural isomorphism
$N^*_m{\cal F}\rightarrow T^*_{f(m)}B$, and,
under this isomorphism, the transversal geodesic flow $G_t$ on
$N^*{\cal F}$ corresponds to the geodesic flow $T^*B$
(see, for instance, \cite{ONeil,Re}).

\end{document}